\documentclass[a4paper,11pt]{amsart}
\usepackage{amscd,amsthm,amsfonts,latexsym,amssymb}
\usepackage[dvips]{graphicx}

\usepackage{bez123,calc,curves,ebezier,epic,eepic,graphicx,multiply,rotating}

\newtheorem{theorem}{Theorem} [section]
\newtheorem{lemma}[theorem]{Lemma}
\newtheorem{corollary}[theorem]{Corollary}

\newtheorem{remark}[theorem]{Remark}

\newcommand{\Ee}{{\mathcal E}}

\newcommand{\ra}{\rightarrow}

\renewcommand{\phi}{\varphi}

\newcommand{\Ga}{\Gamma}
\newcommand{\ve}{\varepsilon}
\newcommand{\si}{\sigma}
\newcommand{\Si}{\Sigma}

\newcommand{\wt}{\widetilde} 
\newcommand{\ov}{\overline} 
\newcommand{\wh}{\widehat} 

\begin{document}

\title{An algorithm to prescribe the configuration of a finite graph}
\author{Paul Baird}
\author{Marius Tiba}
\thanks{The first author thanks the Agence National de Recherche, project ANR-07-BLAN-0251-01, for financial support.} 

\address{Universit\'e europ\'eenne de Bretagne, Universit\'e de Brest; CNRS, Laboratoire de Math\'ematiques--UMR 6205, 6 av.\ Victor Le Gorgeu -- CS 93837, 29238 Brest Cedex, France}

\address{International Computer High School, SOS Mihai Bravu 428, Sector 3, Bucharest 030328, Romania}
\email{Paul.Baird@univ-brest.fr, sirtibamarius@yahoo.com}

\begin{abstract} We provide algorithms involving edge slides, for a connected simple graph to evolve in a finite number of steps to another connected simple graph in a prescribed configuration, and for the regularization of such a graph by the minimization of an appropriate energy functional.  
\end{abstract}



\maketitle


\thispagestyle{empty}

\section{Introduction} \label{intro}
Edge operations on graphs such as edge rotations, or switchings, have been widely used in various contexts: to provide a notion of distance between graphs, for example by measuring the number of moves required to get from one graph to another, see for example \cite{Ch-Sa-Zo, Go-Sw, Ja, Jo}; to provide algorithms to transform one planar graph or tree into another, see \cite{Ai-Re}, and \cite{Bo-Hu} for an overview; as a tool in random graph theory to estimate probabilities, as exploited in \cite{Mc-Wo, Wo}.  We may view such operations as part of a \emph{dynamic} theory, whereby a graph evolves in discrete steps into a different configuration.  In this article we are motivated by geometric considerations to study how a graph may evolve into one that is as regular as possible, or into a prescribed configuration.

In smooth geometry, curvature can be detected by the convergence or divergence of nearby geodesics.  The Gauss-Bonnet Theorem gives an expression for the total curvature of a compact surface in terms of its Euler characteristic.  On a planar graph $\Ga$ there is a notion of \emph{combinatorial curvature} introduced by Y. Higuchi \cite{Hi}, given by the function $\Phi$ defined at each vertex $x$ by the formula:
$$
\Phi (x) = 1 - \frac{d(x)}{2} + \sum_{\si} \frac{1}{|\si |}\,,
$$
where $d(x)$ denotes the degree of $x$, that is the number of edges incident with $x$, and where the sum is taken over all polygons $\si$ incident with $x$, with $|\si |$ representing the number of sides of $\si$.  For a finite connected planar graph, the total curvature is given by $\sum_x\Phi (x) = 2$.  As a variant, on an arbitrary graph, we can take the integer-valued function $K$ defined at each vertex $x$ by $K(x) = 2 - d(x)$.  Up to a multiple of $2$, this may be viewed as an approximation of $\Phi$ on a sparse graph, that is a graph with few connections.   Then we quickly deduce an analogue of the Gauss-Bonnet Theorem: $\sum_{x} K(x) = 2 \chi (\Ga )$, where the sum is taken over all vertices and where $\chi (\Ga )$ is the Euler characteristic given by $n - e$ with $n$ the number of vertices and $e$ the number of edges of $\Ga$.  The problems we address here in a combinatorial setting are akin to the problem of uniformization and that of prescribing the scalar curvature on a surface, as described for example, in \cite{Au}.

An edge slide is an operation whereby we slide one end of an edge along another edge, so a triple of the form $x \sim y \sim z$ with $x \not\sim z$ becomes either $x \sim z \sim y$ or $ y \sim x \sim z$, where the notation $x\sim y$ means that $x$ and $y$ are adjacent vertices.  An edge slide is a purely local operation which preserves the connectedness and the Euler characteristic of a graph.  A basic question is to know whether one can transform one graph into another by edge slides.  Given two connected simple graphs on the same number of vertices and edges, this was shown to be possible by M. Johnson \cite{Jo}.  His proof is non-constructive in that it doesn't provide an algorithm to carry out the required sequence of edge slides.  Our first result provides a constructive proof of a slightly more general theorem.  

Let $\Ga = (V, E)$ and $\Si = (W, F)$ be two finite connected simple graphs with the same number of vertices and edges and let $\psi : V \ra W$ be a bijection between the vertex sets; we will refer to $\Ga$ as the \emph{initial configuration} and $\Si$ as the \emph{prescribed configuration}.  Then we show there is a combination of edges slides on $\Ga$ to produce a new graph $\wt{\Ga} = (V, \wt{E})$ such that $\psi : \wt{\Ga} \ra \Si$ is an isomorphism of graphs.  Furthermore, we provide an algorithm to carry out the sequence of edges slides.

A first step is to show that any single move of an edge which preserves connectedness can be achieved by edge slides, furthermore we provide a specific algorithm to do this.  Our strategy is then to select the vertex $y$ of smallest degree in $\Si$ and to increase (if necessary) the degree of the corresponding vertex $x$ in $\Ga$ until it has degree $n-1$.  This is done by taking a spanning tree and simultaneously evolving both the spanning tree and the graph by edge slides.  We then have to remove edges incident with $x$ in an appropriate way until its degree coincides with that of $y$. Our objective is to remove each of $x$ and $y$ from $\Ga$ and $\Si$, respectively and so to reduce the problem to graphs of successively smaller size.  A difficulty that may arise is that the complements of $x$ and $y$ may not be connected.  This is overcome by first making a judicious choice of moves in both graphs.  

Recall that a degree sequence is a list of non-negative integers $(d_1, \ldots , d_n)$.  The theorem of Erd\H{o}s and Gallai \cite{Er-Ga} gives the conditions when such a sequence has a realization as the degrees of the vertices of a simple graph.  Two non-isomorphic graphs may have the same degree sequence.  A regular graph is one with all vertices of the same degree.  The problem of generating regular graphs, or more generally, graphs of restricted degree sequences, is an important aspect of the theory of random graphs, see \cite{Mc-Wo, Wo}.  If we view the degree sequence as a combinatorial analogue of the metric and take $K$ as above for the curvature, then edge sliding exchanges degrees from one vertex to a neighbouring vertex and as such, appears as an analogue of the evolution of a metric by its curvature.  

We introduce a natural energy functional $\Ee(\Ga )$ associated to a graph $\Ga$ which measures its discrepancy from being regular.  We then describe an algorithm to regularize $\Ga$ by edge slides, which at each step decreases $\Ee$. At the end of the algorithm, the graph is in what we call an \emph{almost regular configuration}, which is as close to being regular as possible given the number of vertices and edges.  

Finally, we allow the creation of vertices and edges in such a way as to preserve the connectivity and the Euler characteristic of a connected simple graph.  We then provide an algorithm for such a graph $\Ga$ to evolve into a prescribed configuration $\Si$ with the same Euler characteristic, where we no longer require $\Ga$ and $\Si$ to have the same number of vertices and edges.

\section{Notation and terminology} \label{sec:slides}  Let $\Ga = (V, E)$ be a simple graph, where $V$ denotes the vertices and $E$ the edges.  Thus we do not allow double edges or loops.  We will write an edge $\ve \in E$ in the form $\ve = \ov{xy}$ when we want to indicate the end points $x,y \in V$, otherwise we also write $x \sim y$ to indicate that $x$ and $y$ are joined by an edge and we will say that $x$ and $y$ are \emph{adjacent} or are \emph{neighbours}.  We set $d(x)$ to be the number of edges incident with the vertex $x$.  Then $\sum_{x \in V} d(x) = 2e$, where $e = |E|$ is the cardinality of $E$.  We now define the operation of sliding.

Let $x,y,z \in V$ be three distinct vertices such that $x \sim y \sim z$ and $x\not\sim z$.  We call such a configuration of vertices a \emph{triple}.  Then a \emph{(simple) slide along the edge $\ov{yz}$ with pivot $x$} creates the new configuration $x\sim z \sim y$.  In most situations, we shall simply say that we slide the edge $\ov{xy}$ to $\ov{xz}$.  Similarly, a slide along the edge $\ov{xy}$ with pivot $z$ creates the new configuration $y \sim x \sim z$.  Both moves are possible since $x \not\sim z$.  In general, we shall refer to the process of performing a sequence of edge slides as \emph{sliding}.

\medskip
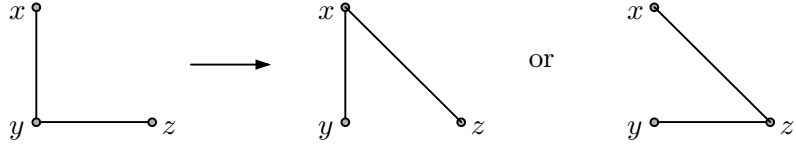
\begin{center}
\setlength{\unitlength}{0.254mm}
\begin{picture}(400,72)(30,-106)
        \allinethickness{0.254mm}\path(40,-40)(40,-100) 
        \allinethickness{0.254mm}\path(40,-100)(100,-100) 
        \allinethickness{0.254mm}\path(200,-40)(200,-100) 
        \allinethickness{0.254mm}\path(360,-100)(420,-100) 
        \allinethickness{0.254mm}\path(120,-70)(160,-70)\special{sh 1}\path(160,-70)(154,-68)(154,-70)(154,-72)(160,-70) 
        \allinethickness{0.254mm}\special{sh 0.3}\put(40,-40){\ellipse{4}{4}} 
        \allinethickness{0.254mm}\special{sh 0.3}\put(40,-100){\ellipse{4}{4}} 
        \allinethickness{0.254mm}\special{sh 0.3}\put(100,-100){\ellipse{4}{4}} 
        \allinethickness{0.254mm}\special{sh 0.3}\put(200,-40){\ellipse{4}{4}} 
        \allinethickness{0.254mm}\special{sh 0.3}\put(200,-100){\ellipse{4}{4}} 
        \allinethickness{0.254mm}\special{sh 0.3}\put(260,-100){\ellipse{4}{4}} 
        \allinethickness{0.254mm}\special{sh 0.3}\put(360,-40){\ellipse{4}{4}} 
        \allinethickness{0.254mm}\special{sh 0.3}\put(360,-100){\ellipse{4}{4}} 
        \allinethickness{0.254mm}\special{sh 0.3}\put(420,-100){\ellipse{4}{4}} 
        \put(26,-46){\shortstack{$x$}} 
        \put(26,-106){\shortstack{$y$}} 
        \put(105,-106){\shortstack{$z$}} 
        \put(186,-46){\shortstack{$x$}} 
        \put(186,-106){\shortstack{$y$}} 
        \put(265,-106){\shortstack{$z$}} 
        \put(346,-46){\shortstack{$x$}} 
        \put(346,-106){\shortstack{$y$}} 
        \put(425,-106){\shortstack{$z$}} 
        \put(295,-71){\shortstack{or}} 
        \allinethickness{0.254mm}\path(200,-40)(260,-100) 
        \allinethickness{0.254mm}\path(360,-40)(420,-100) 
\end{picture}
\end{center}
\medskip

An edge slide always preserves connectedness of a graph and in the case of an edge slide along $\ov{yz}$ with pivot $x$, the degree at $y$ decreases by one and the degree at $z$ increases by one, while all other degrees remain the same.  Note also that sliding is a reversible operation and on allowing the trivial slide which involves no move, we have an equivalence relation on the set of graphs with $n$ vertices and $e$ edges.  More precisely, we shall say that two graphs $\Ga$ and $\Si$ are \emph{slide-equivalent} if there exists a sequence of edge slides in $\Ga$ which yields a new graph $\wt{\Ga}$ such that $\wt{\Ga}$ and $\Si$ are isomorphic, that is, there exists a bijection $\psi : V \ra W$ from the vertices of $\wt{\Ga}$ to those of $\Si$ such that $x \sim y$ if and only if $\psi (x) \sim \psi (y)$.

We write a path joining two vertices $x$ and $u$ in the form $\si = [x:u]$, or if we want to indicate the vertices along the path, by $\si = [x:s_1\cdots s_k:u]$, where $x = s_1$, $u = s_k$ and $s_i \sim s_{i+1}$ for $i = 1, \ldots , k-1$, or finally by $\si = [s_1s_2 \cdots s_k]$ if we do not wish to specify the end points.  

Let $\si : [s_1s_2 \cdots s_k]$ be a path and let $y$ be a vertex not contained in $\si$.  Suppose that $y \sim s_1$ and $y \not\sim s_i$ for $i = 2, \ldots , s_k$.  Then we can perform consecutive slides of the edge $\ov{ys_1}$, first to $\ov{ys_2}$, then to $\ov{ys_3}$ and so on, until we reach $\ov{ys_k}$.  In this case we shall say that we \emph{slide the edge $\ov{ys_1}$ along the path $\si$ to the edge $\ov{ys_k}$}. 

Suppose, given the path $\si$ as above, that $y \sim s_i$ and that $y \not\sim s_j$ for some $i,j \in \{ 1, \ldots , k\}$.  Then whatever other connections exist between $y$ and the vertices of $\si$, it is clear that by a combination of slides we may achieve the configuration whereby $y\not\sim s_i$ and $y \sim s_j$, with all other connections unchanged.  We call such a move a \emph{shuffle of $\ov{ys_i}$ to $\ov{ys_j}$ (along $\si$ with pivot $y$)}.  Shuffling is a useful composite move that we will employ in the next section.  

\medskip
\begin{center}
\setlength{\unitlength}{0.200mm}
\begin{picture}(574,137)(43,-206)
        \allinethickness{0.254mm}\path(45,-140)(75,-170) 
        \allinethickness{0.254mm}\path(75,-170)(115,-190) 
        \allinethickness{0.254mm}\path(115,-190)(165,-200) 
        \allinethickness{0.254mm}\special{sh 0.3}\put(45,-140){\ellipse{4}{4}} 
        \allinethickness{0.254mm}\special{sh 0.3}\put(75,-170){\ellipse{4}{4}} 
        \allinethickness{0.254mm}\special{sh 0.3}\put(115,-190){\ellipse{4}{4}} 
        \allinethickness{0.254mm}\special{sh 0.3}\put(165,-200){\ellipse{4}{4}} 
        \allinethickness{0.254mm}\path(215,-190)(165,-200) 
        \allinethickness{0.254mm}\path(215,-190)(255,-170) 
        \allinethickness{0.254mm}\path(255,-170)(285,-140) 
        \allinethickness{0.254mm}\special{sh 0.3}\put(215,-190){\ellipse{4}{4}} 
        \allinethickness{0.254mm}\special{sh 0.3}\put(255,-170){\ellipse{4}{4}} 
        \allinethickness{0.254mm}\special{sh 0.3}\put(285,-140){\ellipse{4}{4}} 
        \allinethickness{0.254mm}\special{sh 0.3}\put(165,-85){\ellipse{4}{4}} 
        \allinethickness{0.254mm}\path(360,-140)(390,-170) 
        \allinethickness{0.254mm}\path(390,-170)(435,-190) 
        \allinethickness{0.254mm}\path(435,-190)(490,-200) 
        \allinethickness{0.254mm}\path(490,-200)(545,-190) 
        \allinethickness{0.254mm}\path(545,-190)(585,-170) 
        \allinethickness{0.254mm}\path(585,-170)(615,-140) 
        \allinethickness{0.254mm}\special{sh 0.3}\put(360,-140){\ellipse{4}{4}} 
        \allinethickness{0.254mm}\special{sh 0.3}\put(390,-170){\ellipse{4}{4}} 
        \allinethickness{0.254mm}\special{sh 0.3}\put(435,-190){\ellipse{4}{4}} 
        \allinethickness{0.254mm}\special{sh 0.3}\put(490,-200){\ellipse{4}{4}} 
        \allinethickness{0.254mm}\special{sh 0.3}\put(545,-190){\ellipse{4}{4}} 
        \allinethickness{0.254mm}\special{sh 0.3}\put(585,-170){\ellipse{4}{4}} 
        \allinethickness{0.254mm}\special{sh 0.3}\put(615,-140){\ellipse{4}{4}} 
        \allinethickness{0.254mm}\special{sh 0.3}\put(490,-85){\ellipse{4}{4}} 
        \allinethickness{0.254mm}\path(165,-85)(45,-140) 
        \allinethickness{0.254mm}\path(165,-85)(75,-170) 
        \allinethickness{0.254mm}\path(165,-85)(115,-190) 
        \allinethickness{0.254mm}\path(165,-85)(165,-200) 
        \allinethickness{0.254mm}\path(165,-85)(215,-190) 
        \allinethickness{0.254mm}\path(165,-85)(285,-140) 
        \allinethickness{0.254mm}\path(490,-85)(360,-140) 
        \allinethickness{0.254mm}\path(490,-85)(435,-190) 
        \allinethickness{0.254mm}\path(490,-85)(490,-200) 
        \allinethickness{0.254mm}\path(490,-85)(545,-190) 
        \allinethickness{0.254mm}\path(490,-85)(585,-170) 
        \allinethickness{0.254mm}\path(490,-85)(615,-140) 
        \allinethickness{0.254mm}\path(310,-140)(340,-140)\special{sh 1}\path(340,-140)(334,-138)(334,-140)(334,-142)(340,-140) 
        \put(170,-81){\shortstack{$y$}} 
        \put(495,-81){\shortstack{$y$}} 
        \put(65,-186){\shortstack{$s_i$}} 
        \put(260,-186){\shortstack{$s_j$}} 
        \put(380,-186){\shortstack{$s_i$}} 
        \put(590,-181){\shortstack{$s_j$}} 
        \put(215,-206){\shortstack{$s_{j-1}$}} 
        \put(545,-206){\shortstack{$s_{j-1}$}} 
\end{picture}
\end{center}
\medskip

\begin{minipage}[center]{12cm}
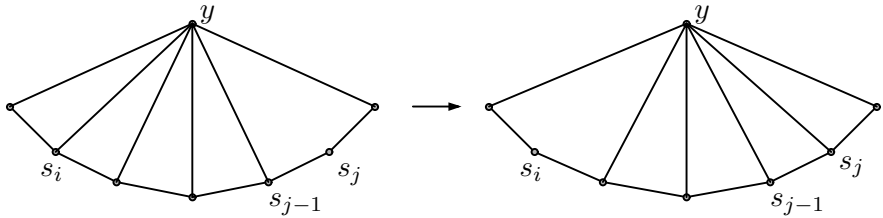
\emph{\footnotesize The figure illustrates a shuffle of the edge $\ov{ys_i}$ to $\ov{ys_j}$; first slide $\ov{ys_{j-1}}$ to $\ov{ys_j}$, then $\ov{ys_{j-2}}$ to $\ov{ys_{j-1}}$ and so on.}
\end{minipage}

\section{Discrete evolution to a prescribed configuration} In this section we show by producing an algorithm, that any two connected simple graphs with the same number of vertices and edges are slide-equivalent.  In fact we establish a slightly stronger result, which allows us to choose arbitrarily the bijection between the vertices.  This is expressed by the following theorem.

\begin{theorem} \label{thm:prescribed}  Given two connected simple graphs $\Ga$ and $\Si$ each with the same number of vertices and edges.  Let $\psi$ be a bijection between the vertices.  Then there exists a combination of edge slides of $\Ga$ to produce a graph $\wt{\Ga}$ such that $\psi : \wt{\Ga} \ra \Si$ is an isomorphism of graphs.
\end{theorem}

As the first part of our strategy to prove this Theorem, we show that any move of an edge which preserves connectivity can be achieved by a combination of slides.  Furthermore, we provide a specific algorithm to do this.  First we have a preliminary lemma.

\begin{lemma} \label{lem:paths}  Let $\Ga = (V, E)$ be a connected simple incomplete graph.  Let $\ve = \ov{uv} \in E$ be an edge and let $x,y \in V$ be unconnected vertices such that the graph obtained from $\Ga$ by removing $\ve$ and adding $\ov{xy}$ is connected.  Then, up to labelling of $u$ and $v$, there exist paths from $x$ to $u$ and from $y$ to $v$ not containing the edge $\ov{uv}$.
\end{lemma}

\begin{proof}  Since $\Ga$ is connected, there is a path from $x$ to the edge $\ve$.  Stop the path at the first vertex of $\ve$ encountered and call this vertex $u$.  This gives a path from $x$ to $u$ not containing $\ve$.  Similarly, $y$ is connected to either $u$ or $v$ by a path not containing $\ve$.  If $y$ is connected to $v$, then we are done.  Suppose on the other hand, $y$ is connected to $u$.  Then by the same argument, since the graph $\wt{\Ga}$ obtained by removing $\ve$ and adding $\ov{xy}$ is connected, there is a path in $\wt{\Ga}$ from $v$ to either $x$ or $y$ not containing $\ov{xy}$.  If this is to $y$, we are done.  If it is to $x$, then on relabeling $u$ and $v$, we are able to obtain paths from $x$ to $u$ and from $y$ to $v$ not containing $\ov{uv}$.
\end{proof}

\begin{lemma} \label{lem:moves}  Any move of an edge which preserves connectedness is a combination of slides.
\end{lemma}

\begin{proof}  Suppose we move $\ve = \ov{uv}$ to $\ov{xy} \ ( x \not\sim y$ in $\Ga$).  Let $\si = [x:s_1\cdots s_k:u]$ be the shortest path from $x$ to $u$ not containing $\ov{uv}$ and let $\tau = [y:t_1 \cdots t_l:v]$ be the shortest path from $y$ to $v$ not containing $\ov{uv}$.  Up to labeling, by Lemma \ref{lem:paths}, these exist.  We have a number of cases to consider.

\medskip

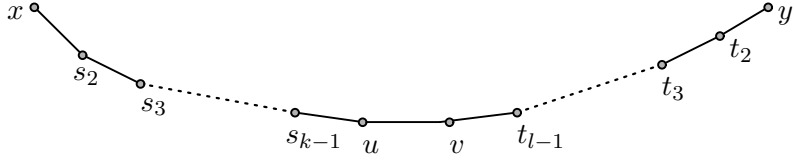
\begin{center}

\setlength{\unitlength}{0.254mm}
\begin{picture}(400,82)(30,-116)
        \allinethickness{0.254mm}\path(40,-40)(65,-65) 
        \allinethickness{0.254mm}\path(65,-65)(95,-80) 
        \allinethickness{0.254mm}\dottedline{5}(95,-80)(175,-95) 
        \allinethickness{0.254mm}\path(175,-95)(210,-100) 
        \allinethickness{0.254mm}\path(210,-100)(250,-100) 
        \allinethickness{0.254mm}\path(250,-100)(290,-95) 
        \allinethickness{0.254mm}\dottedline{5}(290,-95)(365,-70) 
        \allinethickness{0.254mm}\path(365,-70)(395,-55) 
        \allinethickness{0.254mm}\path(395,-55)(420,-40) 
        \allinethickness{0.254mm}\special{sh 0.3}\put(40,-40){\ellipse{4}{4}} 
        \allinethickness{0.254mm}\special{sh 0.3}\put(65,-65){\ellipse{4}{4}} 
        \allinethickness{0.254mm}\special{sh 0.3}\put(95,-80){\ellipse{4}{4}} 
        \allinethickness{0.254mm}\special{sh 0.3}\put(175,-95){\ellipse{4}{4}} 
        \allinethickness{0.254mm}\special{sh 0.3}\put(210,-100){\ellipse{4}{4}} 
        \allinethickness{0.254mm}\special{sh 0.3}\put(255,-100){\ellipse{4}{4}} 
        \allinethickness{0.254mm}\special{sh 0.3}\put(365,-70){\ellipse{4}{4}} 
        \allinethickness{0.254mm}\special{sh 0.3}\put(395,-55){\ellipse{4}{4}} 
        \allinethickness{0.254mm}\special{sh 0.3}\put(420,-40){\ellipse{4}{4}} 
        \allinethickness{0.254mm}\special{sh 0.3}\put(290,-95){\ellipse{4}{4}} 
        \put(26,-46){\shortstack{$x$}} 
        \put(60,-78){\shortstack{$s_2$}} 
        \put(170,-111){\shortstack{$s_{k-1}$}} 
        \put(210,-116){\shortstack{$u$}} 
        \put(255,-116){\shortstack{$v$}} 
        \put(290,-111){\shortstack{$t_{l-1}$}} 
        \put(365,-86){\shortstack{$t_3$}} 
        \put(95,-94){\shortstack{$s_3$}} 
        \put(400,-66){\shortstack{$t_2$}} 
        \put(425,-46){\shortstack{$y$}} 
\end{picture}
\end{center}
\medskip

\noindent \emph{Case 1}: $x = u$; then $\si$ is the trivial path.

(a) $u \not\in\tau$:  then slide or shuffle $\ov{uv}$ to $\ov{xy}$ along $\tau$ (with pivot $u = x$).

\medskip
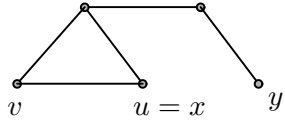
\begin{center}
\setlength{\unitlength}{0.254mm}
\begin{picture}(152,58)(40,-136)
        \allinethickness{0.254mm}\path(80,-80)(140,-80) 
        \allinethickness{0.254mm}\path(110,-120)(80,-80) 
        \allinethickness{0.254mm}\special{sh 0.3}\put(140,-80){\ellipse{4}{4}} 
        \allinethickness{0.254mm}\special{sh 0.3}\put(80,-80){\ellipse{4}{4}} 
        \allinethickness{0.254mm}\special{sh 0.3}\put(110,-120){\ellipse{4}{4}} 
        \allinethickness{0.254mm}\path(80,-80)(45,-120) 
        \allinethickness{0.254mm}\path(140,-80)(170,-120) 
        \allinethickness{0.254mm}\special{sh 0.3}\put(170,-120){\ellipse{4}{4}} 
        \allinethickness{0.254mm}\special{sh 0.3}\put(45,-120){\ellipse{4}{4}} 
        \allinethickness{0.254mm}\path(45,-120)(110,-120) 
        \put(105,-136){\shortstack{$u=x$}} 
        \put(175,-131){\shortstack{$y$}} 
        \put(40,-136){\shortstack{$v$}} 
\end{picture}
\end{center}
\medskip  

(b)  $u \in \tau$: Suppose $u = t_i$; slide $\ov{vu}$ to $\ov{vy}$ along $[u:t_it_{i-1}\cdots t_1:y]$; this is possible since if $v \sim t_j$, $j < i$, then there would be a shorter path from $y$ to $v$.  Now slide $\ov{vy}$ along $[v:t_lt_{l-1}\cdots t_i:x]$ to $\ov{xy}$.  Once more, this is possible, since if there exists $j$, $j>i$, with $y \sim t_j$, then there would be a shorter path from $y$ to $v$.

\medskip
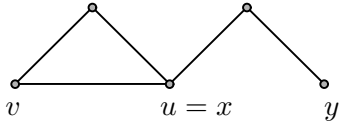
\begin{center}
\setlength{\unitlength}{0.254mm}
\begin{picture}(182,58)(35,-116)
        \allinethickness{0.254mm}\path(40,-100)(80,-60) 
        \allinethickness{0.254mm}\path(80,-60)(120,-100) 
        \allinethickness{0.254mm}\path(120,-100)(160,-60) 
        \allinethickness{0.254mm}\path(160,-60)(200,-100) 
        \allinethickness{0.254mm}\path(40,-100)(120,-100) 
        \allinethickness{0.254mm}\special{sh 0.3}\put(40,-100){\ellipse{4}{4}} 
        \allinethickness{0.254mm}\special{sh 0.3}\put(80,-60){\ellipse{4}{4}} 
        \allinethickness{0.254mm}\special{sh 0.3}\put(120,-100){\ellipse{4}{4}} 
        \allinethickness{0.254mm}\special{sh 0.3}\put(160,-60){\ellipse{4}{4}} 
        \allinethickness{0.254mm}\special{sh 0.3}\put(200,-100){\ellipse{4}{4}} 
        \put(35,-116){\shortstack{$v$}} 
        \put(115,-116){\shortstack{$u=x$}} 
        \put(200,-116){\shortstack{$y$}} 
\end{picture}
\end{center}

\medskip

\noindent \emph{Case 2}:  $x = v$:  

(a)  $u \in \tau$:  Suppose $u = t_i$; then we just slide $\ov{vu}$ to $\ov{vy}$ along $[u:t_it_{i-1}\cdots t_1:y]$ (with pivot $x = v$).

\medskip
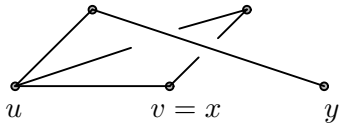
\begin{center}
\setlength{\unitlength}{0.254mm}
\begin{picture}(182,58)(35,-96)
        \allinethickness{0.254mm}\special{sh 0.3}\put(40,-80){\ellipse{4}{4}} 
        \allinethickness{0.254mm}\special{sh 0.3}\put(120,-80){\ellipse{4}{4}} 
        \allinethickness{0.254mm}\special{sh 0.3}\put(200,-80){\ellipse{4}{4}} 
        \allinethickness{0.254mm}\special{sh 0.3}\put(80,-40){\ellipse{4}{4}} 
        \allinethickness{0.254mm}\special{sh 0.3}\put(160,-40){\ellipse{4}{4}} 
        \allinethickness{0.254mm}\path(80,-40)(200,-80) 
        \allinethickness{0.254mm}\path(80,-40)(40,-80) 
        \allinethickness{0.254mm}\path(40,-80)(120,-80) 
        \allinethickness{0.254mm}\path(120,-80)(135,-65) 
        \allinethickness{0.254mm}\path(145,-55)(160,-40) 
        \allinethickness{0.254mm}\path(40,-80)(100,-60) 
        \allinethickness{0.254mm}\path(125,-50)(160,-40) 
        \put(200,-96){\shortstack{$y$}} 
        \put(110,-96){\shortstack{$v=x$}} 
        \put(35,-96){\shortstack{$u$}} 
\end{picture}
\end{center}
\medskip

(b)  $u \not\in \tau$:  If $y = s_i \in \si$, then clearly $x = v \not\in [y:s_is_{i+1}\cdots s_k:u]$, so we can shuffle edge $\ov{vu}$ to $\ov{vy} = \ov{xy}$. 

\medskip
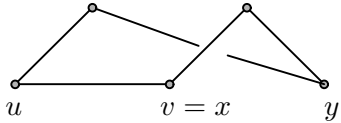
\begin{center}
\setlength{\unitlength}{0.254mm}
\begin{picture}(182,58)(35,-136)
        \allinethickness{0.254mm}\path(40,-120)(80,-80) 
        \allinethickness{0.254mm}\path(40,-120)(120,-120) 
        \allinethickness{0.254mm}\path(120,-120)(160,-80) 
        \allinethickness{0.254mm}\path(160,-80)(200,-120) 
        \allinethickness{0.254mm}\path(80,-80)(135,-100) 
        \allinethickness{0.254mm}\special{sh 0.3}\put(40,-120){\ellipse{4}{4}} 
        \allinethickness{0.254mm}\special{sh 0.3}\put(120,-120){\ellipse{4}{4}} 
        \allinethickness{0.254mm}\special{sh 0.3}\put(80,-80){\ellipse{4}{4}} 
        \allinethickness{0.254mm}\special{sh 0.3}\put(160,-80){\ellipse{4}{4}} 
        \allinethickness{0.254mm}\special{sh 0.3}\put(200,-120){\ellipse{4}{4}} 
        \allinethickness{0.254mm}\path(150,-105)(200,-120) 
        \put(35,-136){\shortstack{$u$}} 
        \put(115,-136){\shortstack{$v=x$}} 
        \put(200,-136){\shortstack{$y$}} 
\end{picture}
\end{center}
\medskip 

If $y\not\in \si$, then we shuffle edge $\ov{uv}$ to $\ov{uy}$ and then $\ov{uy}$ to $\ov{xy}$.

\medskip
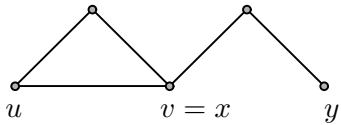
\begin{center}
\setlength{\unitlength}{0.254mm}
\begin{picture}(182,58)(35,-116)
        \allinethickness{0.254mm}\path(40,-100)(80,-60) 
        \allinethickness{0.254mm}\path(80,-60)(120,-100) 
        \allinethickness{0.254mm}\path(120,-100)(160,-60) 
        \allinethickness{0.254mm}\path(160,-60)(200,-100) 
        \allinethickness{0.254mm}\path(40,-100)(120,-100) 
        \allinethickness{0.254mm}\special{sh 0.3}\put(200,-100){\ellipse{4}{4}} 
        \allinethickness{0.254mm}\special{sh 0.3}\put(160,-60){\ellipse{4}{4}} 
        \allinethickness{0.254mm}\special{sh 0.3}\put(120,-100){\ellipse{4}{4}} 
        \allinethickness{0.254mm}\special{sh 0.3}\put(80,-60){\ellipse{4}{4}} 
        \allinethickness{0.254mm}\special{sh 0.3}\put(40,-100){\ellipse{4}{4}} 
        \put(35,-116){\shortstack{$u$}} 
        \put(115,-116){\shortstack{$v=x$}} 
        \put(200,-116){\shortstack{$y$}} 
\end{picture}
\end{center}
\medskip

By symmetry, the cases $y = v$ and $y = u$ are dealt with similarly.

\medskip

\noindent \emph{Case 3}:  $x \neq u,v$, $y \neq u,v$: If $x\sim v$, we first slide $\ov{xv}$ to $\ov{xy}$;  then from Case 1, we can now slide $\ov{uv}$ to $\ov{xv}$.  If $x\not\sim v$, then we first slide $\ov{uv}$ to $\ov{xv}$ and once more, from Case 1, we can slide $\ov{xv}$ to $\ov{xy}$.   
\end{proof}

The algorithm to perform a move of an edge by slides is therefore as follows:  \emph{first construct the shortest paths from $x$ to $u$ and from $y$ to $v$ not containing $\ov{uv}$ (interchanging the labels $u$ and $v$ if necessary); then, depending on the case that occurs, proceed as in the above proof.} 

By interchanging two vertices in a graph, we mean the following.  Let $\Ga$ be a graph and let $x,y$ be two vertices of $\Ga$, then \emph{interchanging $x$ and $y$} gives the new graph $\wt{\Ga}$ with all neighbours of $x$ now neighbours of $y$ and all neighbours of $y$ now neighbours of $x$; all other edges remain unchanged. The following lemma shows that interchanging vertices can be achieved by sliding.

\begin{lemma}  \label{lem:interchange}  In any connected simple graph, interchanging two vertices $x$ and $y$ can be achieved by a combination of slides.
\end{lemma}

\begin{proof}  \emph{Case 1}:  $x \sim y$:  If $z$ is a vertex such that $ z \sim x$ and $z \sim y$, then there is no slide needed.  If $z \sim x$ and $z\not\sim y$, then we slide $\ov{zx}$ to $\ov{zy}$.  

\medskip

\noindent \emph{Case 2}:  $ x \not\sim y$:  Let $\si = [x:s_1\cdots s_l:y]$ be the shortest path joining $x$ to $y$; clearly $l\geq 3$.  If $z\not\in \si$ is a vertex such that $z \sim x$ and $z\not\sim y$, then we can shuffle $\ov{zx}$ to $\ov{zy}$.  Similarly, if $z\not\in \si$ with $z\not\sim x$ and $z\sim y$, we can shuffle $\ov{zy}$ to $\ov{zx}$.  There remain the edges $\ov{xs_2}$ and $\ov{s_{l-1}y}$ to deal with.

We observe that since $\si$ has minimal length, for any $3\leq i \leq l$, we have $x \not\sim s_i\not\sim y$.  If $l = 3$, there is nothing left to do.  If $l \geq 4$, then by sliding, we have to move edge $\ov{xs_2}$ to $\ov{ys_2}$ and $\ov{s_{l-1}y}$ to $\ov{s_{l-1}x}$ to finish the problem.  We do this as in the diagram below:  

\medskip

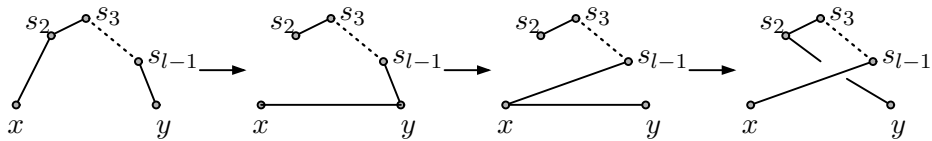
\begin{center}

\setlength{\unitlength}{0.230mm}
\begin{picture}(542,77)(30,-136)
        \allinethickness{0.254mm}\path(35,-120)(55,-80) 
        \allinethickness{0.254mm}\path(55,-80)(75,-70) 
        \allinethickness{0.254mm}\path(115,-120)(105,-95) 
        \allinethickness{0.254mm}\dottedline{5}(75,-70)(105,-95) 
        \allinethickness{0.254mm}\path(140,-100)(165,-100)\special{sh 1}\path(165,-100)(159,-98)(159,-100)(159,-102)(165,-100) 
        \allinethickness{0.254mm}\path(195,-80)(215,-70) 
        \allinethickness{0.254mm}\path(255,-120)(245,-95) 
        \allinethickness{0.254mm}\dottedline{5}(215,-70)(245,-95) 
        \allinethickness{0.254mm}\special{sh 0.3}\put(35,-120){\ellipse{4}{4}} 
        \allinethickness{0.254mm}\special{sh 0.3}\put(55,-80){\ellipse{4}{4}} 
        \allinethickness{0.254mm}\special{sh 0.3}\put(75,-70){\ellipse{4}{4}} 
        \allinethickness{0.254mm}\special{sh 0.3}\put(105,-95){\ellipse{4}{4}} 
        \allinethickness{0.254mm}\special{sh 0.3}\put(115,-120){\ellipse{4}{4}} 
        \allinethickness{0.254mm}\special{sh 0.3}\put(175,-120){\ellipse{4}{4}} 
        \allinethickness{0.254mm}\special{sh 0.3}\put(195,-80){\ellipse{4}{4}} 
        \allinethickness{0.254mm}\special{sh 0.3}\put(215,-70){\ellipse{4}{4}} 
        \allinethickness{0.254mm}\special{sh 0.3}\put(245,-95){\ellipse{4}{4}} 
        \allinethickness{0.254mm}\special{sh 0.3}\put(255,-120){\ellipse{4}{4}} 
        \allinethickness{0.254mm}\path(175,-120)(255,-120) 
        \allinethickness{0.254mm}\path(280,-100)(305,-100)\special{sh 1}\path(305,-100)(299,-98)(299,-100)(299,-102)(305,-100) 
        \allinethickness{0.254mm}\path(315,-120)(395,-120) 
        \allinethickness{0.254mm}\path(385,-95)(315,-120) 
        \allinethickness{0.254mm}\path(335,-80)(355,-70) 
        \allinethickness{0.254mm}\dottedline{5}(355,-70)(385,-95) 
        \allinethickness{0.254mm}\path(525,-95)(455,-120) 
        \allinethickness{0.254mm}\path(510,-105)(535,-120) 
        \allinethickness{0.254mm}\path(475,-80)(495,-95) 
        \allinethickness{0.254mm}\path(475,-80)(495,-70) 
        \allinethickness{0.254mm}\dottedline{5}(495,-70)(525,-95) 
        \allinethickness{0.254mm}\path(420,-100)(445,-100)\special{sh 1}\path(445,-100)(439,-98)(439,-100)(439,-102)(445,-100) 
        \allinethickness{0.254mm}\special{sh 0.3}\put(315,-120){\ellipse{4}{4}} 
        \allinethickness{0.254mm}\special{sh 0.3}\put(385,-95){\ellipse{4}{4}} 
        \allinethickness{0.254mm}\special{sh 0.3}\put(335,-80){\ellipse{4}{4}} 
        \allinethickness{0.254mm}\special{sh 0.3}\put(355,-70){\ellipse{4}{4}} 
        \allinethickness{0.254mm}\special{sh 0.3}\put(395,-120){\ellipse{4}{4}} 
        \allinethickness{0.254mm}\special{sh 0.3}\put(455,-120){\ellipse{4}{4}} 
        \allinethickness{0.254mm}\special{sh 0.3}\put(475,-80){\ellipse{4}{4}} 
        \allinethickness{0.254mm}\special{sh 0.3}\put(495,-70){\ellipse{4}{4}} 
        \allinethickness{0.254mm}\special{sh 0.3}\put(525,-95){\ellipse{4}{4}} 
        \allinethickness{0.254mm}\special{sh 0.3}\put(535,-120){\ellipse{4}{4}} 
        \put(30,-136){\shortstack{$x$}} 
        \put(115,-136){\shortstack{$y$}} 
        \put(80,-71){\shortstack{$s_3$}} 
        \put(110,-96){\shortstack{$s_{l-1}$}} 
        \put(170,-136){\shortstack{$x$}} 
        \put(255,-136){\shortstack{$y$}} 
        \put(220,-71){\shortstack{$s_3$}} 
        \put(250,-96){\shortstack{$s_{l-1}$}} 
        \put(310,-136){\shortstack{$x$}} 
        \put(395,-136){\shortstack{$y$}} 
        \put(360,-71){\shortstack{$s_3$}} 
        \put(390,-96){\shortstack{$s_{l-1}$}} 
        \put(450,-136){\shortstack{$x$}} 
        \put(535,-136){\shortstack{$y$}} 
        \put(530,-96){\shortstack{$s_{l-1}$}} 
        \put(40,-76){\shortstack{$s_2$}} 
        \put(180,-76){\shortstack{$s_2$}} 
        \put(320,-76){\shortstack{$s_2$}} 
        \put(460,-76){\shortstack{$s_2$}} 
        \put(500,-71){\shortstack{$s_3$}} 
\end{picture}

\end{center}
\medskip

First, slide $\ov{xs_2}$ to $\ov{xy}$; now slide $\ov{s_{l-1}y}$ to $\ov{s_{l-1}x}$, and finally, slide $\ov{yx}$ to $\ov{ys_2}$.  
\end{proof} 

\begin{lemma} \label{lem:spanning}  Let $T$ be a tree on $n$ vertices (and so $n-1$ edges).  Fix any vertex, then there exists a combination of slides which give this vertex the degree $n-1$.
\end{lemma}

\begin{proof}  Note that since sliding does not change the number of vertices or edges of a graph, it must preserve the property that a graph be a tree.  

Clearly the result is true for $n = 1,2$, so suppose $n \geq 3$.  Fix a vertex $x$.  If $d(x) = n-1$ we are done, so suppose $d(x)<n-1$.  There must be at least two vertices of degree $1$, otherwise we would contradict the fact that the sum of degrees is equal to $2(n-1)$.  Furthermore, there must be at least one such vertex which is not a neighbour of $x$.  Let $y$ be a vertex of degree $1$ not equal to $x$ and not a neighbour of $x$.  Let $z$ be its neighbour.  Then we can shuffle $\ov{yz}$ to $\ov{yx}$ so increasing the degree at $x$ by one.  Repeat this process until $d(x) = n-1$.
\end{proof}

\begin{lemma} \label{lem:increasing-degree}  Let $\Ga$ be a connected simple graph and let $x$ be any chosen vertex.  Then by sliding we can arrange for $x$ to have degree $n-1$.
\end{lemma}

\begin{proof} Consider the graph $\Ga_1 = \Ga$.  Fix a vertex $x$.  If $d(x) = n-1$ we are done, otherwise, take a spanning tree $T_1$ of $\Ga_1$.  Now begin the operation of the proof of Lemma \ref{lem:spanning} by performing a simple slide of an edge of $T_1$ in order to increase the degree at $x$. If the graph permits (so as not to create a double edge), slide the same edge in $\Ga_1$, otherwise do nothing.  We now have a new tree $T_2$ spanning a new graph $\Ga_2$ (maybe $\Ga_1 = \Ga_2$).  Repeat this process, so that at each step we obtain a connected simple graph $\Ga_i$ obtained from $\Ga_{i-1}$ by an edge-slide, together with a spanning tree $T_i$, in such a way that $T_i$ is obtained from $T_{i-1}$ by performing the requisite simple slide to increase the degree of the designated vertex.  We continue until the designated vertex has degree $n - 1$.     
\end{proof}

\medskip 

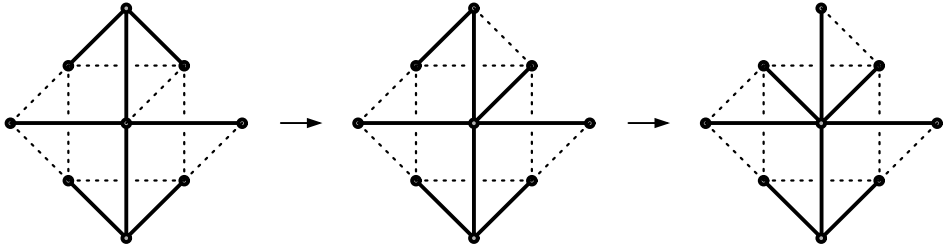
\setlength{\unitlength}{0.254mm}
\begin{picture}(484,124)(38,-162)
        \allinethickness{0.508mm}\path(40,-100)(160,-100) 
        \allinethickness{0.508mm}\path(100,-40)(100,-160) 
        \allinethickness{0.508mm}\path(100,-40)(130,-70) 
        \allinethickness{0.508mm}\path(100,-40)(70,-70) 
        \allinethickness{0.508mm}\path(100,-160)(130,-130) 
        \allinethickness{0.508mm}\path(100,-160)(70,-130) 
        \allinethickness{0.508mm}\special{sh 0.3}\put(160,-100){\ellipse{4}{4}} 
        \allinethickness{0.508mm}\special{sh 0.3}\put(130,-70){\ellipse{4}{4}} 
        \allinethickness{0.508mm}\special{sh 0.3}\put(100,-40){\ellipse{4}{4}} 
        \allinethickness{0.508mm}\special{sh 0.3}\put(70,-70){\ellipse{4}{4}} 
        \allinethickness{0.508mm}\special{sh 0.3}\put(40,-100){\ellipse{4}{4}} 
        \allinethickness{0.508mm}\special{sh 0.3}\put(70,-130){\ellipse{4}{4}} 
        \allinethickness{0.508mm}\special{sh 0.3}\put(100,-160){\ellipse{4}{4}} 
        \allinethickness{0.508mm}\special{sh 0.3}\put(130,-130){\ellipse{4}{4}} 
        \allinethickness{0.508mm}\special{sh 0.3}\put(100,-100){\ellipse{4}{4}} 
        \allinethickness{0.254mm}\dottedline{5}(160,-100)(130,-130) 
        \allinethickness{0.254mm}\dottedline{5}(70,-130)(40,-100) 
        \allinethickness{0.254mm}\dottedline{5}(40,-100)(70,-70) 
        \allinethickness{0.254mm}\dottedline{5}(70,-70)(95,-70) 
        \allinethickness{0.254mm}\dottedline{5}(105,-70)(130,-70) 
        \allinethickness{0.254mm}\dottedline{5}(70,-70)(70,-95) 
        \allinethickness{0.254mm}\dottedline{5}(70,-105)(70,-130) 
        \allinethickness{0.254mm}\dottedline{5}(70,-130)(95,-130) 
        \allinethickness{0.254mm}\dottedline{5}(105,-130)(130,-130) 
        \allinethickness{0.254mm}\dottedline{5}(130,-130)(130,-105) 
        \allinethickness{0.254mm}\dottedline{5}(130,-95)(130,-70) 
        \allinethickness{0.254mm}\dottedline{5}(130,-70)(100,-100) 
        \allinethickness{0.254mm}\path(180,-100)(200,-100)\special{sh 1}\path(200,-100)(194,-98)(194,-100)(194,-102)(200,-100) 
        \allinethickness{0.508mm}\path(220,-100)(340,-100) 
        \allinethickness{0.508mm}\path(280,-40)(280,-160) 
        \allinethickness{0.508mm}\path(280,-40)(250,-70) 
        \allinethickness{0.508mm}\path(280,-160)(310,-130) 
        \allinethickness{0.508mm}\path(280,-160)(250,-130) 
        \allinethickness{0.508mm}\path(280,-100)(310,-70) 
        \allinethickness{0.508mm}\special{sh 0.3}\put(340,-100){\ellipse{4}{4}} 
        \allinethickness{0.508mm}\special{sh 0.3}\put(310,-70){\ellipse{4}{4}} 
        \allinethickness{0.508mm}\special{sh 0.3}\put(280,-40){\ellipse{4}{4}} 
        \allinethickness{0.508mm}\special{sh 0.3}\put(280,-100){\ellipse{4}{4}} 
        \allinethickness{0.508mm}\special{sh 0.3}\put(250,-70){\ellipse{4}{4}} 
        \allinethickness{0.508mm}\special{sh 0.3}\put(220,-100){\ellipse{4}{4}} 
        \allinethickness{0.508mm}\special{sh 0.3}\put(250,-130){\ellipse{4}{4}} 
        \allinethickness{0.508mm}\special{sh 0.3}\put(280,-160){\ellipse{4}{4}} 
        \allinethickness{0.508mm}\special{sh 0.3}\put(310,-130){\ellipse{4}{4}} 
        \allinethickness{0.254mm}\dottedline{5}(250,-70)(220,-100) 
        \allinethickness{0.254mm}\dottedline{5}(280,-40)(310,-70) 
        \allinethickness{0.254mm}\dottedline{5}(340,-100)(310,-130) 
        \allinethickness{0.254mm}\dottedline{5}(250,-130)(220,-100) 
        \allinethickness{0.254mm}\dottedline{5}(250,-130)(250,-105) 
        \allinethickness{0.254mm}\dottedline{5}(250,-95)(250,-70) 
        \allinethickness{0.254mm}\dottedline{5}(250,-70)(275,-70) 
        \allinethickness{0.254mm}\dottedline{5}(285,-70)(310,-70) 
        \allinethickness{0.254mm}\dottedline{5}(310,-70)(310,-95) 
        \allinethickness{0.254mm}\dottedline{5}(310,-105)(310,-130) 
        \allinethickness{0.254mm}\dottedline{5}(310,-130)(285,-130) 
        \allinethickness{0.254mm}\dottedline{5}(275,-130)(250,-130) 
        \allinethickness{0.254mm}\path(360,-100)(380,-100)\special{sh 1}\path(380,-100)(374,-98)(374,-100)(374,-102)(380,-100) 
        \allinethickness{0.508mm}\path(400,-100)(520,-100) 
        \allinethickness{0.508mm}\path(460,-40)(460,-160) 
        \allinethickness{0.508mm}\path(460,-100)(490,-70) 
        \allinethickness{0.508mm}\path(460,-100)(430,-70) 
        \allinethickness{0.508mm}\path(460,-160)(490,-130) 
        \allinethickness{0.508mm}\path(460,-160)(430,-130) 
        \allinethickness{0.508mm}\special{sh 0.3}\put(460,-160){\ellipse{4}{4}} 
        \allinethickness{0.508mm}\special{sh 0.3}\put(490,-130){\ellipse{4}{4}} 
        \allinethickness{0.508mm}\special{sh 0.3}\put(430,-130){\ellipse{4}{4}} 
        \allinethickness{0.508mm}\special{sh 0.3}\put(460,-100){\ellipse{4}{4}} 
        \allinethickness{0.508mm}\special{sh 0.3}\put(400,-100){\ellipse{4}{4}} 
        \allinethickness{0.508mm}\special{sh 0.3}\put(520,-100){\ellipse{4}{4}} 
        \allinethickness{0.508mm}\special{sh 0.3}\put(490,-70){\ellipse{4}{4}} 
        \allinethickness{0.508mm}\special{sh 0.3}\put(430,-70){\ellipse{4}{4}} 
        \allinethickness{0.508mm}\special{sh 0.3}\put(460,-40){\ellipse{4}{4}} 
        \allinethickness{0.254mm}\dottedline{5}(460,-40)(490,-70) 
        \allinethickness{0.254mm}\dottedline{5}(430,-70)(400,-100) 
        \allinethickness{0.254mm}\dottedline{5}(400,-100)(430,-130) 
        \allinethickness{0.254mm}\dottedline{5}(490,-130)(520,-100) 
        \allinethickness{0.254mm}\dottedline{5}(430,-70)(455,-70) 
        \allinethickness{0.254mm}\dottedline{5}(465,-70)(490,-70) 
        \allinethickness{0.254mm}\dottedline{5}(430,-70)(430,-95) 
        \allinethickness{0.254mm}\dottedline{5}(430,-105)(430,-130) 
        \allinethickness{0.254mm}\dottedline{5}(430,-130)(455,-130) 
        \allinethickness{0.254mm}\dottedline{5}(465,-130)(490,-130) 
        \allinethickness{0.254mm}\dottedline{5}(490,-130)(490,-105) 
        \allinethickness{0.254mm}\dottedline{5}(490,-95)(490,-70) 
\end{picture}

\medskip
\begin{minipage}[center]{12cm}
\emph{\footnotesize The figure illustrates the simultaneous evolution of a graph (bold and dotted lines) and its spanning tree (bold lines).  At the first step, the top right hand edge of the spanning tree slides to increase the degree of the central vertex by one, while the graph remains unchanged.  At the second step, both the graph and the tree evolve to increase the degree of the central vertex by one.}
\end{minipage}
\bigskip 

\noindent \emph{Proof of Theorem} \ref{thm:prescribed}.  We will apply induction on the number $n$ of vertices in $\Ga$.

\medskip

\noindent \emph{Step 1}:  In $\Si$, consider a vertex $y$ of minimum degree $d_1$.  Let $x = \psi^{-1}(y)$ be the corresponding vertex in $\Ga$.  By Lemma \ref{lem:increasing-degree}, by sliding, we can obtain a new graph $\wt{\Ga}$ in which $x$ has degree $n - 1$.  Let $C_1, \ldots , C_s$ be the connected components of the graph formed from $\wt{\Ga}$ be deleting $x$ and all its incident edges.  By connectivity of $\wt{\Ga}$, each one of these components is connected to $x$ by one or more edges.  Let $k_i$ be the number of vertices of $C_i$.  Our aim is to decrease the degree of $x$ until it reaches $d_1$.  To do this, we perform the following algorithm.

\medskip

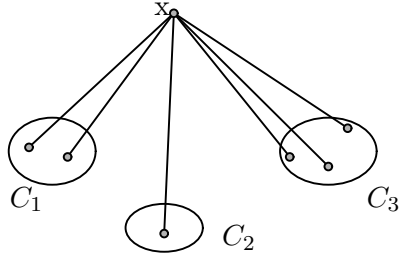
\begin{center}

\setlength{\unitlength}{0.254mm}
\begin{picture}(216,136)(80,-145)
        \allinethickness{0.254mm}\put(102,-92){\ellipse{45}{35}} 
        \allinethickness{0.254mm}\put(160,-132){\ellipse{40}{25}} 
        \allinethickness{0.254mm}\put(245,-92){\ellipse{50}{35}} 
        \allinethickness{0.254mm}\path(165,-20)(90,-90) 
        \allinethickness{0.254mm}\path(165,-20)(110,-95) 
        \allinethickness{0.254mm}\path(165,-20)(160,-135) 
        \allinethickness{0.254mm}\path(165,-20)(225,-95) 
        \allinethickness{0.254mm}\path(165,-20)(245,-100) 
        \allinethickness{0.254mm}\path(165,-20)(255,-80) 
        \allinethickness{0.254mm}\special{sh 0.3}\put(90,-90){\ellipse{4}{4}} 
        \allinethickness{0.254mm}\special{sh 0.3}\put(110,-95){\ellipse{4}{4}} 
        \allinethickness{0.254mm}\special{sh 0.3}\put(165,-20){\ellipse{4}{4}} 
        \allinethickness{0.254mm}\special{sh 0.3}\put(160,-135){\ellipse{4}{4}} 
        \allinethickness{0.254mm}\special{sh 0.3}\put(225,-95){\ellipse{4}{4}} 
        \allinethickness{0.254mm}\special{sh 0.3}\put(245,-100){\ellipse{4}{4}} 
        \allinethickness{0.254mm}\special{sh 0.3}\put(255,-80){\ellipse{4}{4}} 
        \put(155,-21){\shortstack{x}} 
        \put(80,-121){\shortstack{$C_1$}} 
        \put(190,-141){\shortstack{$C_2$}} 
        \put(265,-121){\shortstack{$C_3$}} 
\end{picture}
\end{center}
\medskip

If $s \geq 2$, take an edge from $x$ to $C_1$ and use it to join $C_1$ to $C_2$.  This maintains connectivity and by Lemma \ref{lem:moves}, this can be achieved by a number of slides.  Now $C_1$ and $C_2$ have become one connected component.  We repeat this process until \emph{either} we achieve the degree $d_1$, \emph{or} $d(x)>d_1$ and there exists just one connected component $C$ on removing $x$ and its incident edges.  We claim that in the latter case there exists at least two non-adjacent vertices in $C$.

If $C$ is complete, it has $n - 1$ vertices and $(n-1)(n-2)/2$ edges.  Now the number of edges in $\wt{\Ga}$ is $\geq d_1 + 1 + \frac{(n-1)(n-2)}{2}$.  But, by removing $y$ from $\Si$ together with its incident edges, we see that $\Si$ has at most $d_1 + \frac{(n-1)(n-2)}{2}$ edges.  But this is a contradiction, since all of our operations are achieved by sliding, which preserves the number of edges.  Thus there exist two non-adjacent vertices in $C$ to which we can move an edge from $x$ to $C$.  We now repeat this process until the degree $d_1$ is achieved at the vertex $x$.

We now apply Lemma \ref{lem:interchange}.  By interchanging vertices, we can arrange that the set of vertices that are neighbours of $x$ by the algorithm of Step 1, are exactly those vertices that are required to be neighbours in order that $x\sim y \Rightarrow \psi (x) \sim \psi (y)$.  
Write the new graph that results from Step 1 as $\wh{\Ga}$. 

\medskip

\noindent \emph{Step 2}: Let $\Ga^{\prime}$ be the graph obtained from $\wh{\Ga}$ by removing $x$ and its incident edges.  We aim to apply induction, but $\Ga^{\prime}$ may not be connected.  Let $L_1, \ldots , L_t$ be its connected components.  We will now connect them by slides.

Let $L_i$ have $\ell_i$ vertices.  If $L_1$ contains $\geq \ell_1$ edges, it contains a cycle and we can take one of its edges and join it in $\wh{\Ga}$ to another connected component, $L_2$ say, without disconnecting $L_1$.  We repeat this process.  If finally, we have components $L_1, \ldots , L_r$ ($r\geq 2$) in $\Ga^{\prime}$ each with $\ell_1, \ldots , \ell_r$ vertices and $\ell_1-1, \ldots , \ell_r - 1$ edges (where, for simplicity of notation, we maintain the symbols $\Ga^{\prime}$ and $L_i$, even though these may have changed under the above operations), then the total number of edges is given by $(\ell_1-1)+ \cdots + (\ell_r - 1) + d_1 = n - 1 - r + d_1$.  On the other hand, $|E| \geq d_1n/2$, since $d_1$ is the smallest degree in $\Si$ and $\Si$ and $\wh{\Ga}$ have the same number of edges.  Thus
$$
d_1 \leq \frac{2(n-1-r)}{n-2}<2\,,
$$
so that $d_1 = 1$.  But this is a contradiction, since $r \leq d_1$.  Thus $r = 1$ and $\Ga^{\prime} = L_1$ is connected.

Similarly, in the graph $\Si$, we perform the above operations in order that, on removing $y$ and its incident edges, the resulting graph $\Si^{\prime}$ is connected.  Note that the operations of Step 2 have no effect on the neighbours of $x$ and $y$; these are preserved.  As remarked in Section \ref{sec:slides}, sliding is a reversible operation, hence the slides performed in $\Si$ can be reversed.  By induction, the theorem is true for the graphs $\Ga^{\prime}$ and $\Si^{\prime}$ on $n-1$ vertices and the result follows.
\hfill  $\Box$     

\section{Regularisation} \label{sec:reg} 
Given two positive integers $n$ and $e$, then by Euclidean division, there are unique $k$ and $r$ such that $2e = nk + r$ where $0 \leq r \leq n-1$.  Then we have the following configuration of degrees:
\begin{equation} \label{min-config}
\begin{array}{rl}
 r & \left\{  \begin{array}{l} k+1 \\ \vdots \\ k+1 \end{array} \right.  \\  
n - r & \left\{  \begin{array}{l} k \\ \vdots \\ k \end{array} \right.
\end{array}
\end{equation} 
whose sum is $2e$.  We can characterize such a configuration in the following way.  For a degree sequence $(d_1, \ldots , d_n)$, define its \emph{energy} to be the quantity $\Ee = \sum_i d_i{}^2$.   If each $d_i>0$, we shall say that the sequence is \emph{positive}. 

\begin{lemma}  Given $e>0$, a positive degree sequence $(d_1, \ldots , d_n)$ with sum $2e$ minimizes $\Ee$ amongst all other positive sequences with sum $2e$ if and only if it has the form {\rm (\ref{min-config})} for some $k\geq 1$ and $r$ with $0 \leq r \leq n-1$.  
\end{lemma} 
\begin{proof}  Among all sums $d_1{}^2 + \cdots + d_n{}^2$ with $d_1 + \cdots + d_n = 2e$ $(d_i >0)$, there is one which is minimal; let this be $S$.  Let the sequence $(d_1, \ldots , d_n)$ realize $S$.  We want to show this has the form (\ref{min-config}).  If there are two integers $d_i$ and $d_j$ such that $d_i - d_j \geq 2$, then  the sequence $(d_1, \ldots d_{j - 1}, d_j +1, d_{j+1}, \ldots d_{i-1}, d_i-1, d_{i+1}, \ldots d_n)$ has all terms $>0$, sum $2e$, but with sum of squares equal to $S + 2(d_j - d_i + 1) < S-1$, contradicting our hypothesis that $S$ is the minimum of $\Ee$ amongst positive sequences.  This means that any two terms of the sequence $(d_1, \ldots , d_n)$ can either be equal or can differ by $1$, so they all have the form $t$ or $t+1$ for some positive integer $t$.  But then the sum of the terms lies between $nt$ and $n(t+1)$ and $t$ must equal $k$.    
\end{proof}

Call a graph $\Ga = (V, E)$ such that $\max_{x,y \in V} |d(x)- d(y)| \leq 1$ an \emph{almost regular graph}.  We now devise an algorithm for obtaining an almost regular graph from a given connected simple graph by sliding.   

\medskip

\begin{center}

\setlength{\unitlength}{0.254mm}
\begin{picture}(404,84)(38,-102)
        \allinethickness{0.254mm}\path(40,-60)(120,-60) 
        \allinethickness{0.254mm}\path(80,-20)(80,-100) 
        \allinethickness{0.254mm}\path(80,-20)(120,-60) 
        \allinethickness{0.254mm}\path(120,-60)(80,-100) 
        \allinethickness{0.254mm}\path(80,-100)(40,-60) 
        \allinethickness{0.254mm}\path(40,-60)(80,-20) 
        \allinethickness{0.254mm}\path(100,-40)(60,-80) 
        \allinethickness{0.254mm}\path(100,-80)(60,-40) 
        \allinethickness{0.254mm}\path(200,-60)(280,-60) 
        \allinethickness{0.254mm}\path(240,-20)(240,-100) 
        \allinethickness{0.254mm}\path(240,-20)(280,-60) 
        \allinethickness{0.254mm}\path(280,-60)(240,-100) 
        \allinethickness{0.254mm}\path(240,-100)(200,-60) 
        \allinethickness{0.254mm}\path(200,-60)(240,-20) 
        \allinethickness{0.254mm}\path(220,-40)(260,-80) 
        \allinethickness{0.254mm}\path(260,-40)(245,-40) 
        \allinethickness{0.254mm}\path(235,-40)(220,-40) 
        \allinethickness{0.254mm}\path(260,-80)(245,-80) 
        \allinethickness{0.254mm}\path(235,-80)(220,-80) 
        \allinethickness{0.254mm}\path(400,-20)(400,-100) 
        \allinethickness{0.254mm}\path(400,-20)(440,-60) 
        \allinethickness{0.254mm}\path(400,-20)(360,-60) 
        \allinethickness{0.254mm}\path(360,-60)(400,-100) 
        \allinethickness{0.254mm}\path(400,-100)(440,-60) 
        \allinethickness{0.254mm}\path(380,-40)(400,-100) 
        \allinethickness{0.254mm}\path(420,-80)(400,-20) 
        \allinethickness{0.254mm}\path(360,-60)(380,-60) 
        \allinethickness{0.254mm}\path(390,-60)(395,-60) 
        \allinethickness{0.254mm}\path(405,-60)(410,-60) 
        \allinethickness{0.254mm}\path(420,-60)(440,-60) 
        \allinethickness{0.254mm}\path(380,-40)(395,-40) 
        \allinethickness{0.254mm}\path(410,-40)(420,-40) 
        \allinethickness{0.254mm}\path(420,-80)(405,-80) 
        \allinethickness{0.254mm}\path(390,-80)(380,-80) 
        \allinethickness{0.254mm}\path(150,-60)(170,-60)\special{sh 1}\path(170,-60)(164,-58)(164,-60)(164,-62)(170,-60) 
        \allinethickness{0.254mm}\path(310,-60)(330,-60)\special{sh 1}\path(330,-60)(324,-58)(324,-60)(324,-62)(330,-60) 
        \allinethickness{0.254mm}\special{sh 0.3}\put(120,-60){\ellipse{4}{4}} 
        \allinethickness{0.254mm}\special{sh 0.3}\put(100,-40){\ellipse{4}{4}} 
        \allinethickness{0.254mm}\special{sh 0.3}\put(80,-20){\ellipse{4}{4}} 
        \allinethickness{0.254mm}\special{sh 0.3}\put(60,-40){\ellipse{4}{4}} 
        \allinethickness{0.254mm}\special{sh 0.3}\put(40,-60){\ellipse{4}{4}} 
        \allinethickness{0.254mm}\special{sh 0.3}\put(60,-80){\ellipse{4}{4}} 
        \allinethickness{0.254mm}\special{sh 0.3}\put(80,-100){\ellipse{4}{4}} 
        \allinethickness{0.254mm}\special{sh 0.3}\put(100,-80){\ellipse{4}{4}} 
        \allinethickness{0.254mm}\special{sh 0.3}\put(80,-60){\ellipse{4}{4}} 
        \allinethickness{0.254mm}\special{sh 0.3}\put(240,-20){\ellipse{4}{4}} 
        \allinethickness{0.254mm}\special{sh 0.3}\put(220,-40){\ellipse{4}{4}} 
        \allinethickness{0.254mm}\special{sh 0.3}\put(200,-60){\ellipse{4}{4}} 
        \allinethickness{0.254mm}\special{sh 0.3}\put(220,-80){\ellipse{4}{4}} 
        \allinethickness{0.254mm}\special{sh 0.3}\put(240,-100){\ellipse{4}{4}} 
        \allinethickness{0.254mm}\special{sh 0.3}\put(260,-80){\ellipse{4}{4}} 
        \allinethickness{0.254mm}\special{sh 0.3}\put(280,-60){\ellipse{4}{4}} 
        \allinethickness{0.254mm}\special{sh 0.3}\put(260,-40){\ellipse{4}{4}} 
        \allinethickness{0.254mm}\special{sh 0.3}\put(240,-60){\ellipse{4}{4}} 
        \allinethickness{0.254mm}\special{sh 0.3}\put(400,-20){\ellipse{4}{4}} 
        \allinethickness{0.254mm}\special{sh 0.3}\put(380,-40){\ellipse{4}{4}} 
        \allinethickness{0.254mm}\special{sh 0.3}\put(360,-60){\ellipse{4}{4}} 
        \allinethickness{0.254mm}\special{sh 0.3}\put(380,-80){\ellipse{4}{4}} 
        \allinethickness{0.254mm}\special{sh 0.3}\put(400,-100){\ellipse{4}{4}} 
        \allinethickness{0.254mm}\special{sh 0.3}\put(420,-80){\ellipse{4}{4}} 
        \allinethickness{0.254mm}\special{sh 0.3}\put(440,-60){\ellipse{4}{4}} 
        \allinethickness{0.254mm}\special{sh 0.3}\put(420,-40){\ellipse{4}{4}} 
        \allinethickness{0.254mm}\special{sh 0.3}\put(400,-60){\ellipse{4}{4}} 
\end{picture}

\end{center}
\medskip

\begin{minipage}[center]{12cm}
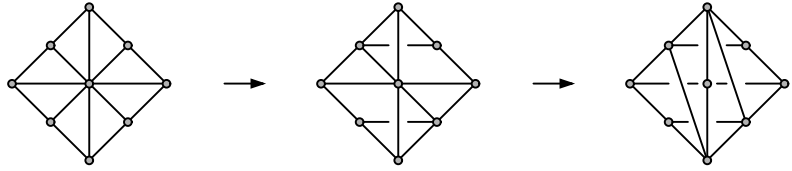
\emph{\footnotesize The figure shows the regularization according to the algorithm below of a graph on nine vertices and sixteen edges.  Sliding occurs in such a way that degrees are exchanged from high incidence vertices to low incidence ones until an almost regular configuration is achieved.  At each step two slides have occured.}
\end{minipage}
\bigskip 

\begin{theorem}  Let $\Ga = (V, E)$ be a connected simple graph.  Then $\Ga$ is slide-equivalent to an almost regular graph.
\end{theorem}

\begin{proof}  If no two vertices have degrees which differ by at least $2$, then the graph is already in an almost regular configuration, so we suppose otherwise.  We make the following algorithm:  At each step, we take two vertices $x,y \in V$ such that $d(x) - d(y) \geq 2$.  Let $\si = [x:s_1s_2\cdots s_k:y]$ be the shortest path from $x$ to $y$.

If $k = 2$, then $x \sim y$ and because $d(x) - d(y) \geq 2$, there exists a vertex $a$ such that $a\sim x$,  $a \not\sim y$ and $a \neq y$.  Then we slide the edge $\ov{ax}$ along $\ov{xy}$, so decreasing the degree at $x$ by $1$, increasing the degree at $y$ by $1$ and leaving all other degrees the same.

If $k \geq 3$, then $x\not\sim y$ and because $d(x) - d(y) \geq 2$, there are two distinct vertices $a,b$ such that $x\sim a$, $x \sim b$, $a \neq y\neq b$ and $a\not\sim y\not\sim b$.  Since $\si$ is the shortest path from $x$ to $y$, we cannot have both $a$ and $b$ in this path; so suppose that $a \not\in \si$.  Then we shuffle the edge $\ov{ax}$ along $\si$ until we connect $a$ to $y$.

Now we see that after each step, the energy $\Ee$ strictly decreases and so this algorithm ends.  In fact, for a connected simple graph, the degrees satisfy $1\leq d(x) \leq n-1$ at each vertex $x$, so there is only a finite number of possible degree sequences and consequently, only a finite number of possible values of $\Ee$.  The graph we obtain when we finish has vertices all of whose degrees differ by at most $1$ and so is almost regular, as required.
\end{proof} 

\begin{remark}  The above theorem shows that any degree sequence of the form (\ref{min-config}) where $n-1\leq e\leq n(n-1)/2$ can be realized as the degree sequence of a connected simple graph.  We simply construct any connected simple graph with $e$ edges and then apply the theorem to slide it into an almost regular configuration.
\end{remark}  

\section{Expanding and collapsing a graph within its Euler class} \label{sec:exp}
Theorem \ref{thm:prescribed} shows that any two connected simple graphs with the same number of edges and vertices are equivalent by edge slides.  However, it would be useful to be able to change the number of vertices and edges and still establish equivalence under appropriate conditions.  This may be important in random graph theory, where one is required to let $n \ra\infty$ in an appropriate class of graphs, for example regular graphs, see \cite{Bo, Mc-Wo, Wo}.  In this section, we show how it is possible to increase the number of vertices and edges without limit whilst preserving the Euler characteristic.

Let $\Ga = (V, E)$ be a finite graph.  Recall that the Euler characteristic is the number $\chi (\Ga ) = n - e$, where $n$ is the number of vertices and $e$ the number of edges.  Note that if $\Ga$ is connected and simple (no loops or double edges), then there are two extremes; a complete graph or a tree.  For a complete graph $e = \frac{1}{2}n(n-1)$ and for a tree $e = n-1$, so that 
$$
n-1 \leq e \leq \frac{1}{2}n(n-1)
$$
and we deduce that
$$
\frac{n(3-n)}{2} \leq \chi \leq 1\,.
$$

We now introduce a process for expanding and collapsing a graph.  This involves two possible moves which are equivalent up to sliding.  Furthermore, the moves preserve the Euler characteristic and the connectivity of the graph.

\medskip

\noindent \emph{Expanding}:

(i)  We add a new vertex $y$ to the graph which is attached to any given vertex $x$ by a new edge $\overline{xy}$ (so that $d(y) = 1$). 

(ii)  To any edge $\overline{xz}$, we add a new vertex $y$ to its interior, so dividing the edge into two edges $\overline{zy}$ and $\overline{yx}$.

\medskip

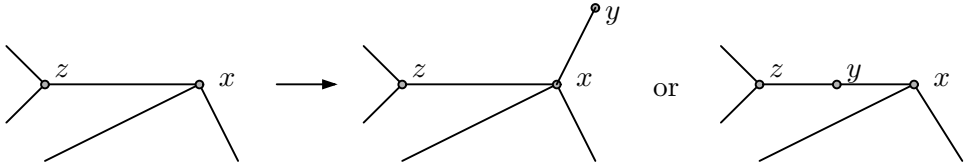
\begin{center}
\setlength{\unitlength}{0.254mm}
\begin{picture}(497,86)(60,-140)
        \allinethickness{0.254mm}\path(80,-100)(160,-100) 
        \allinethickness{0.254mm}\path(60,-80)(80,-100) 
        \allinethickness{0.254mm}\path(80,-100)(60,-120) 
        \allinethickness{0.254mm}\path(160,-100)(80,-140) 
        \allinethickness{0.254mm}\path(160,-100)(180,-140) 
        \allinethickness{0.254mm}\path(265,-100)(345,-100) 
        \allinethickness{0.254mm}\path(265,-100)(245,-80) 
        \allinethickness{0.254mm}\path(265,-100)(245,-120) 
        \allinethickness{0.254mm}\path(345,-100)(265,-140) 
        \allinethickness{0.254mm}\path(345,-100)(365,-140) 
        \allinethickness{0.254mm}\path(450,-100)(530,-100) 
        \allinethickness{0.254mm}\path(450,-100)(430,-80) 
        \allinethickness{0.254mm}\path(450,-100)(430,-120) 
        \allinethickness{0.254mm}\path(530,-100)(450,-140) 
        \allinethickness{0.254mm}\path(530,-100)(555,-140) 
        \allinethickness{0.254mm}\special{sh 0.3}\put(80,-100){\ellipse{4}{4}} 
        \allinethickness{0.254mm}\special{sh 0.3}\put(160,-100){\ellipse{4}{4}} 
        \allinethickness{0.254mm}\special{sh 0.3}\put(265,-100){\ellipse{4}{4}} 
        \allinethickness{0.254mm}\special{sh 0.3}\put(345,-100){\ellipse{4}{4}} 
        \allinethickness{0.254mm}\special{sh 0.3}\put(365,-60){\ellipse{4}{4}} 
        \allinethickness{0.254mm}\special{sh 0.3}\put(450,-100){\ellipse{4}{4}} 
        \allinethickness{0.254mm}\special{sh 0.3}\put(490,-100){\ellipse{4}{4}} 
        \allinethickness{0.254mm}\special{sh 0.3}\put(530,-100){\ellipse{4}{4}} 
        \put(85,-96){\shortstack{$z$}} 
        \put(170,-101){\shortstack{$x$}} 
        \put(270,-96){\shortstack{$z$}} 
        \put(355,-101){\shortstack{$x$}} 
        \put(370,-66){\shortstack{$y$}} 
        \put(455,-96){\shortstack{$z$}} 
        \put(495,-96){\shortstack{$y$}} 
        \put(540,-101){\shortstack{$x$}} 
        \put(395,-106){\shortstack{or}} 
        \allinethickness{0.254mm}\path(200,-100)(230,-100)\special{sh 1}\path(230,-100)(224,-98)(224,-100)(224,-102)(230,-100) 
        \allinethickness{0.254mm}\path(365,-60)(345,-100) 
\end{picture}
\end{center}

\medskip

Both these operations preserve the Euler characteristic, since we add one edge and one vertex.  Furthermore, if $z \sim x$, then after (i), we can slide $\overline{zx}$ to $\overline{zy}$, so that effectively, we have just added a vertex to the interior of the edge $\overline{zx}$.  Conversely, if we add a vertex $y$ to the interior of $\overline{zx}$, then we can slide $\overline{zy}$ to $\overline{zx}$, leaving the isolated vertex $y$ such that $y \sim x$.  Thus (i) and (ii) are equivalent by edge slides. 

\medskip

\noindent \emph{Collapsing}:  

(iii)  If $y$ is a vertex such that $d(y) = 1$, then delete $y$ and the edge that joins it to the rest of the graph.

(iv)  Let $y$ be a vertex such that $d(y) = 2$ and such that if $x\sim y$ and $z \sim y$ ($z \neq x$), then $x\not\sim z$; then delete the vertex $y$ together with the edges $\overline{zy}$ and $\overline{yx}$ and add the edge $\overline{zx}$.

\medskip

Both operations (iii) and (iv) remove one vertex and one edge and so leave the Euler characteristic the same.

\begin{corollary} {\rm (of Theorem \ref{thm:prescribed})}  Let $\Ga$ and $\Si$ be two connected simple graph with $m$ and $n$ vertices, respectively, having the same Euler characteristic.  Suppose $m \leq n$.  Then we can expand $\Ga$ using operations {\rm (i)} or {\rm (ii)} and perform a combination of edge slides to obtain a graph isomorphic to $\Si$.  Equally, we can collapse $\Si$ using operations {\rm (iii)} or {\rm (iv)} and perform a combination of slides to obtain a graph isomorphic to $\Ga$.
\end{corollary}

\begin{proof}  The first part of the corollary is a direct consequence of Theorem \ref{thm:prescribed}:  we add $(n-m)$ vertices and edges according to the rules (i) and (ii) (it doesn't matter where we place them); then we can slide the graph into a configuration that is isomorphic to $\Si$.

For the second part, since the Euler characteristics of $\Ga$ and $\Si$ coincide, we must have
$$
\frac{m(3-m)}{2} \leq \chi \leq 1\,.
$$ 
By induction, it suffices to show that the result holds for $m = n-1$.  Substituting this into the above inequality, we have
\begin{equation} \label{edges}
\frac{(n-1)(4 - n)}{2} \leq \chi \leq 1 \Rightarrow f \leq \frac{(n-1)(n-2)}{2} + 1\,,
\end{equation}
where $f$ is the number of edges of $\Si$.  

Take any vertex $x$ of $\Si$ and let $C$ be the graph obtained from $\Si$ by removing $x$ and its incident edges.  Let $C_1, \ldots , C_r$ be the connected components of $C$.  Note that $r\leq d(x)$.  If $d(x) >1$, then we take an edge joining $x$ to one of its neighbours and move it to two unconnected vertices in $C$.  In the first instance we move such an edge to connect two disconnected components of $C$, if such exist.  Note that there must always be two unconnected vertices in $C$, since if $C$ is complete, then it must contain $(n-1)(n-2)/2$ edges, contradicting the inequality (\ref{edges}).  We repeat this process until $d(x) = 1$.  Let $\widetilde{C}$ be the graph obtained by removing $x$ and the one edge incident with it.  Then $\widetilde{C}$ must be connected (since $r \leq d(x)$).  We now collapse the graph $\Si$ by removing the vertex $x$ of degree $1$ and the one edge incident with it, to obtain the graph $\widetilde{C}$ on $n-1$ vertices, with the same Euler characteristic $\chi$.  By Theorem \ref{thm:prescribed}, by sliding in $\widetilde{C}$, we can obtain a graph on $n-1$ vertices isomorphic to $\Ga$.
\end{proof}

\end{document}